\documentclass[a4paper,11pt]{amsart}

\usepackage{graphicx}
\usepackage{mathptmx}
\usepackage{amsmath}
\usepackage{amssymb}
\usepackage{enumitem}
\usepackage{xcolor} 
\usepackage{pgfplots}

\newmuskip\pFqmuskip

\newcommand*\pFq[6][8]{%
  \begingroup 
  \pFqmuskip=#1mu\relax
  \mathcode`=\string"8000
  \begingroup\lccode`\~=`\,
  \lowercase{\endgroup\let~}\pFqcomma
  F^{#2}_{#3}{\left(\genfrac..{0pt}{}{#4}{#5}\bigg|#6\right)}%
  \endgroup
}
\newcommand{\pFqcomma}{\mskip\pFqmuskip}

\newtheorem{theorem}{Theorem}[section]

\begin{document}

\title[]{Spivey-type recurrence relation for fully degenerate Bell polynomials}

\author{Taekyun  Kim}
\address{Department of Mathematics, Kwangwoon University, Seoul 139-701, Republic of Korea}
\email{tkkim@kw.ac.kr}
\author{Dae San  Kim}
\address{Department of Mathematics, Sogang University, Seoul 121-742, Republic of Korea}
\email{dskim@sogang.ac.kr}


\begin{abstract}
Spivey's combinatorial method revealed an important identity for Bell numbers, involving Stirling numbers of the second kind. This paper extends his work by deriving Spivey-type recurrence relations for fully degenerate Bell polynomials and degenerate Fubini polynomials. Our derivation uses degenerate Stirling numbers of the second kind and two-variable degenerate Fubini polynomials of order $\alpha$. 
\end{abstract}

\maketitle

\markboth{\centerline{\scriptsize Spivey-type recurrence relation for fully degenerate Bell polynomials}}
{\centerline{\scriptsize Taekyun Kim and Dae San Kim}}

\section{Introduction} 
Spivey derived the following remarkable identity for Bell numbers using a combinatorial method (see \eqref{8}, [22]):
\begin{equation}
\phi_{n+m}=\sum_{k=0}^{m}\sum_{l=0}^{n}{m \brace k}\binom{n}{l}k^{n-l}\phi_{l},\quad (m,n\ge 0), \label{1}
\end{equation}
which has since inspired extensive research on Spivey-type relations (see [4,6,8,12-17] and the references therein]). \par

Spivey's Bell number formula in \eqref{1} provides new ways to calculate Bell numbers and their generalizations by relating them to Stirling numbers of the second kind through a combinatorial identity. Its implications include the derivation of new sum formulas for Bell numbers, extensions of the formula to more general mathematical objects like $r$-Whitney and degenerate Bell numbers, and the development of combinatorial proofs and interpretations for these generalizations. The following are the key implications of Spivey's formula. \par
\vspace{0.05in}
{\bf{New Sum Formulas:}} 
The formula provides a fresh way to express Bell numbers as sums involving Stirling numbers of the second kind, which leads to new sum formulas for the Bell numbers themselves (see [22]). \par
\vspace{0.05in}
{\bf{Generalizations:}}
The underlying principles of Spivey's formula have been extended to various other combinatorial quantities, including: \\
\indent $\bullet$ $r$-Whitney Numbers: The formula applies to $r$-Whitney numbers, leading to new identities and a deeper understanding of their structure (see [16]). \\
\indent $\bullet$ Degenerate Bell Numbers: Spivey-type recurrences have been established for degenerate Bell and Dowling polynomials, providing new insights into these related sequences (see [12,15,16]). \\
\indent $\bullet$ Lah-Bell Polynomials: Spivey-type realtions have been extended to $r$-Lah-Bell and $\lambda$-analogue of $r$-Lah-Bell polynomials (see [8]). \\
\indent $\bullet$ $q$-Generalizations: The formula has been generalized to the $q$-analogues of Bell and Stirling numbers, opening up connections to advanced topics in combinatorics and $q$-calculus (see [6]). \par
\vspace{0.05in}
{\bf{New Proof Techniques:}} 
The formula has inspired different approaches to proofs, including: \\
\indent $\bullet$ Generating Function Proofs: Researchers have provided generating function proofs for Spivey's result, offering an algebraic perspective (see [8,12-14]). \\
\indent $\bullet$ Use of Boson and Differential Operators: The use of boson and differential operators provides a powerful algebraic framework for deriving Spivey-type relation. This method bypasses traditional combinatorial proofs by leveraging the correspondence between these operators (see [6,8,15]).\\
\indent $\bullet$ Rook Polynomials: The formula has found a connection with rook polynomials, providing novel, bijective proofs for certain combinatorial identities (see [4]). \par
\vspace{0.05in}
{\bf{Combinatorial Interpretations:}}
The formula provides a basis for developing new combinatorial interpretations, helping to explain the meaning of these number sequences in terms of set partitions and other structures (see [22]). \par
\vspace{0.05in}
{\bf{Connections to Other Fields:}}
Its implications extend to other areas, such as: \\
\indent $\bullet$ Probability: Bell numbers can be interpreted as moments of a Poisson distribution, and Spivey's formula contributes to this understanding (see [13,14,17]). \\
\indent $\bullet$ Graphs: The formula can be applied to graph theory for counting partitions where blocks are independent sets (see [4]). \par
\vspace{0.1in}
Recently, Kim-Kim showed the following recurrence relation for the degenerate Bell polynomials (see \eqref{7}), which is given by 
\begin{equation}
\phi_{n+m,\lambda}(x)=\sum_{k=0}^{m}\sum_{l=0}^{n}{m \brace k}_{\lambda}\binom{n}{l}x^{k}\big(k-m\lambda\big)_{n-l,\lambda}\phi_{l,\lambda}(x),\quad (\mathrm{see}\ [11]),\label{2}
\end{equation}
where $m,n$ are nonnegative integers. \\ 
Letting $\lambda \rightarrow 0$ in \eqref{2} gives Spivey's recurrence relation for Bell polynomials (see \eqref{8}):
\begin{equation}
\phi_{n+m}(x)=\sum_{k=0}^{m}\sum_{l=0}^{n}{m \brace k}\binom{n}{l}x^{k}k^{n-l}\phi_{l}(x),\quad (m,n\ge 0). \label{3}
\end{equation}
Moreover, letting $x=1$ yields the Spivey's relation in \eqref{1}. \par
In this paper, we show the following two Spivey-type relations:
\begin{align}
&\mathrm{Bel}_{n+m,\lambda}(x)= \sum_{k=0}^{m}\sum_{l=0}^{n}(1)_{k,\lambda}{m \brace k}_{\lambda}\binom{n}{l}x^{k}F_{n-l,\lambda}^{(k)}(-\lambda x,k-m\lambda)\mathrm{Bel}_{l,\lambda}(x), \label{4} \\
&F_{n+m,\lambda}(x)=\sum_{k=0}^{m}\sum_{l=0}^{n}k!{m \brace k}_{\lambda}\binom{n}{l}x^{k}F_{n-l,\lambda}^{(k)}(x,k-m\lambda)F_{l,\lambda}(x), \nonumber
\end{align}
where $\mathrm{Bel}_{n,\lambda}(x)$ are the fully degenerate Bell polynomials (see \eqref{9}), $F_{n,\lambda}^{(\alpha)}(x,y)$ are the two variable degenerate Fubini polynomials of order $\alpha$ (see \eqref{11}), and $F_{n,\lambda}(x)$ are the degenerate Fubini polynomials (see \eqref{14}).
Here we note that, by letting $\lambda \rightarrow 0$, we get \eqref{3} from \eqref{4} (see \eqref{10}, \eqref{12}).

For any nonzero $\lambda\in\mathbb{R}$, the degenerate exponentials are defined by 
\begin{equation}
e_{\lambda}^{x}(t)=\sum_{k=0}^{\infty}(x)_{k,\lambda}\frac{t^{k}}{k!}=(1+\lambda t)^{\frac{x}{\lambda}},\quad e_{\lambda}(t)=e_{\lambda}^{1}(t),\quad (\mathrm{see}\ [5,9-12,14]), \label{5}
\end{equation}
where 
\begin{displaymath}
(x)_{0,\lambda}=1,\quad (x)_{n,\lambda}=x(x-\lambda)\cdots\big(x-(n-1)\lambda\big),\quad (n\ge 1).  
\end{displaymath}
The degenerate Stirling numbers ${n \brace k}_{\lambda}$ of the second kind are given by (see\ [12])
\begin{align}
&(x)_{n,\lambda}=\sum_{k=0}^{n}{n \brace k}_{\lambda}(x)_{k},\quad (n\ge 0), \label{6} \\
&\frac{1}{k!}\big(e_{\lambda}(t)-1\big)^{k}=\sum_{n=k}^{\infty}{n \brace k}_{\lambda}\frac{t^{n}}{n!},\quad (k\ge 0). \nonumber
\end{align}
Note that 
\begin{displaymath}
\lim_{\lambda\rightarrow 0}{n \brace k}_{\lambda}={n \brace k}, 
\end{displaymath}
where ${n \brace k}$ is the ordinary Stirling number of the second kind (see [3,18]). \par 
The degenerate Bell polynomials are defined by (see \eqref{6}, \ [11,12])
\begin{align}
&\phi_{n,\lambda}(x)=\sum_{k=0}^{n}{n \brace k}_{\lambda}x{^k},\quad (n\ge 0), \label{7} \\
&e^{x(e_{\lambda}(t)-1)}=\sum_{n=0}^{\infty}\phi_{n,\lambda}(x)\frac{t^{n}}{n!}. \nonumber
\end{align}
When $x=1,\ \phi_{n,\lambda}=\phi_{n,\lambda}(1)$ are called the degenerate Bell numbers. 
Note that $\lim_{\lambda\rightarrow 0}\phi_{n,\lambda}(x)=\phi_{n}(x)$, where $\phi_{n}(x)$ are the classical Bell polynomials given by
\begin{align}
&\phi_{n}(x)=\sum_{k=0}^{n}{n \brace k} x{^k},\quad (n\ge 0), \label{8} \\
&e^{x(e^{t}-1)}=\sum_{n=0}^{\infty}\phi_{n}(x)\frac{t^{n}}{n!}. \nonumber
\end{align}
The fully degenerate Bell polynomials are introduced by Kim-Kim as (see \eqref{6}, \ [5])
\begin{align}
&\mathrm{Bel}_{n,\lambda}(x)=\sum_{k=0}^{n}{n \brace k}_{\lambda}(1)_{k,\lambda}x^{k},\quad (n\ge 0), \label{9} \\
&e_{\lambda}\Big(x\big(e_{\lambda}(t)-1\big)\Big)=\sum_{n=0}^{\infty}\mathrm{Bel}_{n,\lambda}(x)\frac{t^{n}}{n!}. \nonumber
\end{align}
We note here that 
\begin{equation} \label{10}
\lim_{\lambda\rightarrow 0} \mathrm{Bel}_{n,\lambda}(x)=\phi_{n}(x),\quad (n\ge 0). 
\end{equation}
When $x=1$, $\mathrm{Bel}_{n,\lambda}= \mathrm{Bel}_{n,\lambda}(1),\ (n\ge 0)$, are called the fully degenerate Bell numbers. \par
Two variable Fubini polynomials of order $\alpha$, $F_{n}^{(\alpha)}(x,y)$, are given by 
\begin{equation*}
\bigg(\frac{1}{1-x(e^{t}-1)}\bigg)^{\alpha}e^{yt}=\sum_{n=0}^{\infty}F_{n}^{(\alpha)}(x,y)\frac{t^{n}}{n!}. 
\end{equation*}
As a degenerate version of $F_{n}^{(\alpha)}(x,y)$, the two variable degenerate Fubini polynomials of order $\alpha $ are introduced by Kim-Kim as
\begin{equation}
\bigg(\frac{1}{1-x(e_{\lambda}(t)-1)}\bigg)^{\alpha}e_{\lambda}^{y}(t)=\sum_{n=0}^{\infty}F_{n,\lambda}^{(\alpha)}(x,y)\frac{t^{n}}{n!},\quad (\mathrm{see}\ [9]). \label{11}
\end{equation}
When $x=0$, we note that 
\begin{equation}
F_{n,\lambda}^{(\alpha)}(0,y)=(y)_{n,\lambda},\quad (n\ge 0). \label{12}	
\end{equation}
By \eqref{11}, we easily get (see \eqref{6})
\begin{equation*}
F_{n,\lambda}^{(\alpha)}(x,0)=\sum_{k=0}^{n}\langle \alpha\rangle_{k}{n \brace k}_{\lambda}x{^k},\quad (n\ge 0), 
\end{equation*}
where 
\begin{displaymath}
\langle x\rangle_{0}=1,\quad \langle x\rangle_{n}=x(x+1)(x+2)\cdots (x+n-1),\ (n\ge 1). 
\end{displaymath}
Note that 
\begin{displaymath}
\lim_{\lambda\rightarrow 0}F_{n,\lambda}^{(1)}(x,0)=\sum_{k=0}^{n}k!{n \brace k}x^{k}=F_{n}(x), 
\end{displaymath}
where $F_{n}(x)$ are the ordinary Fubini polynomials given by (see [11])
\begin{equation}
\frac{1}{1-x(e^{t}-1)}=\sum_{n=0}^{\infty}F_{n}(x)\frac{t^{n}}{n!}.\label{13}
\end{equation} 
The degenerate Fubini polynomials are given by (see \eqref{6}, \ [11])
\begin{align}
&F_{n,\lambda}(x)=\sum_{k=0}^{n}k!{n \brace k}_{\lambda}x^{k},\quad (n\ge 0), \label{14} \\ 
&\frac{1}{1-x(e_{\lambda}(t)-1)}=\sum_{n=0}^{\infty}F_{n,\lambda}(x)\frac{t^{n}}{n!}. \nonumber
\end{align}
Note that 
\begin{displaymath}
\lim_{\lambda\rightarrow 0}F_{n,\lambda}(x)=\sum_{k=0}^{n}k!{n \brace k}x^{k}=F_{n}(x).
\end{displaymath} \par
The reader may refer to [1,2,5,9-12,14] for the recent developments on various degenerate versions of many special numbers and polynomials, and [3,18,20] as general references for this paper.

\section{Spivey-type recurrence relation for fully degenerate Bell polynomials} 
By \eqref{9}, we get 
\begin{align}
&\sum_{m,n=0}^{\infty}\mathrm{Bel}_{n+m,\lambda}\frac{x^{n}}{n!}\frac{y^{m}}{m!}=\sum_{m=0}^{\infty}\bigg(\frac{d}{dx}\bigg)^{m}\sum_{n=0}^{\infty} \mathrm{Bel}_{n,\lambda}\frac{x^{n}}{n!}\frac{y^{m}}{m!}\label{15} \\
&=\sum_{n=0}^{\infty}\frac{\mathrm{Bel}_{n,\lambda}}{n!}\sum_{m=0}^{\infty}\bigg(\frac{d}{dx}\bigg)^{m}x^{n}\frac{y^{m}}{m!}=\sum_{n=0}^{\infty}\frac{\mathrm{Bel}_{n,\lambda}}{n!}\sum_{m=0}^{n}\binom{n}{m}x^{n-m}y^{m} \nonumber\\
&=\sum_{n=0}^{\infty}\frac{\mathrm{Bel}_{n,\lambda}}{n!}(x+y)^{n}=e_{\lambda}\big(e_{\lambda}(x+y)-1\big)=e_{\lambda}\bigg(e_{\lambda}(x)e_{\lambda}\bigg(\frac{y}{1+\lambda x}\bigg)-1\bigg).\nonumber	
\end{align} \par
Now, we observe that 
\begin{align}
&e_{\lambda}\bigg(e_{\lambda}(x)e_{\lambda}\bigg(\frac{y}{1+\lambda x}\bigg)-1\bigg)=\bigg(1+\lambda\bigg(e_{\lambda}(x)e_{\lambda}\bigg(\frac{y}{1+\lambda x}\bigg)-1\bigg)\bigg)^{\frac{1}{\lambda}}\label{16} \\
&=\bigg(1+\lambda e_{\lambda}(x)-\lambda e_{\lambda}(x)+\lambda	 e_{\lambda}(x)e_{\lambda}\bigg(\frac{y}{1+\lambda x}\bigg)-\lambda\bigg)^{\frac{1}{\lambda}}\nonumber\\
&=\bigg(1+\lambda\big(e_{\lambda}(x)-1\big)+\lambda e_{\lambda}(x)\bigg(e_{\lambda}\bigg(\frac{y}{1+\lambda x}\bigg)-1\bigg)\bigg)^{\frac{1}{\lambda}} \nonumber \\
&=\big(1+\lambda(e_{\lambda}(x)-1)\big)^{\frac{1}{\lambda}}\bigg(1+\lambda\frac{e_{\lambda}(x)\big(e_{\lambda}\big(\frac{y}{1+\lambda x}\big)-1\big)}{1+\lambda(e_{\lambda}(x)-1)}\bigg)^{\frac{1}{\lambda}}\nonumber \\
&=e_{\lambda}\big(e_{\lambda}(x)-1\big)e_{\lambda}\bigg(\frac{e_{\lambda}(x)\big(e_{\lambda}\big(\frac{y}{1+\lambda x}\big)-1\big)}{1+\lambda(e_{\lambda}(x)-1)}\bigg)\nonumber\\
&=e_{\lambda}\big(e_{\lambda}(x)-1\big)\sum_{k=0}^{\infty}\frac{(1)_{k,\lambda}}{k!}\bigg(\frac{e_{\lambda}(x)}{1+\lambda(e_{\lambda}(x)-1)}\bigg)^{k}\bigg(e_{\lambda}\bigg(\frac{y}{1+\lambda x}\bigg)-1\bigg)^{k}\nonumber\\
&=e_{\lambda}\big(e_{\lambda}(x)-1\big)\sum_{k=0}^{\infty}(1)_{k,\lambda}\bigg(\frac{1}{1+\lambda(e_{\lambda}(x)-1)}\bigg)^{k}e_{\lambda}^{k}(x)\sum_{m=k}^{\infty}{m \brace k}_{\lambda}e_{\lambda}^{-m\lambda}(x)\frac{y^{m}}{m!}\nonumber\\
&=\sum_{m=0}^{\infty}\frac{y^{m}}{m!}\bigg(\sum_{k=0}^{m}(1)_{k,\lambda}{m \brace k}_{\lambda} \bigg(\frac{1}{1+\lambda(e_{\lambda}(x)-1)}\bigg)^{k}e^{k-m\lambda}(x)e_{\lambda}\big(e_{\lambda}(x)-1\big)\nonumber\\
&=\sum_{m=0}^{\infty}\frac{y^{m}}{m!}\sum_{j=0}^{\infty}\frac{x^{j}}{j!}\sum_{k=0}^{m}(1)_{k,\lambda}{m \brace k}_{\lambda}F_{j,\lambda}^{(k)}(-\lambda,k-m\lambda)\sum_{l=0}^{\infty}\mathrm{Bel}_{l,\lambda}\frac{x^{l}}{l!}\nonumber\\
&=\sum_{m=0}^{\infty}\sum_{n=0}^{\infty}\bigg(\sum_{l=0}^{n}\sum_{k=0}^{m}(1)_{k,\lambda}{m \brace k}_{\lambda}\binom{n}{l}\mathrm{Bel}_{l,\lambda}F_{n-l,\lambda}^{(k)}(-\lambda,k-m \lambda) \bigg)\frac{y^{m}}{m!}\frac{x^{n}}{n!}. \nonumber
\end{align}
Therefore, by \eqref{15} and \eqref{16}, we obtain the following theorem. 
\begin{theorem}
For $m,n\ge 0$, we have 
\begin{displaymath}
\mathrm{Bel}_{n+m,\lambda}=\sum_{l=0}^{n}\sum_{k=0}^{m}(1)_{k,\lambda}{m \brace k}_{\lambda}\binom{n}{l}\mathrm{Bel}_{l,\lambda}F_{n-l,\lambda}^{(k)}(-\lambda,k-m\lambda). 
\end{displaymath}
\end{theorem}
Letting $\lambda \rightarrow 0$ gives Spivey's recurrence relation for Bell numbers (see \eqref{10}, \eqref{12}):
\begin{align*}
\phi_{n+m}&=\lim_{\lambda\rightarrow 0}\mathrm{Bel}_{n+m,\lambda}=\sum_{l=0}^{n}\sum_{k=0}^{m}{m \brace k}\binom{n}{l}\phi_{l}k^{n-l}. 
\end{align*}
From \eqref{9}, we note that 
\begin{align}
&\sum_{n=0}^{\infty}\mathrm{Bel}_{n+m,\lambda}(t)\frac{x^{n}}{n!}\frac{y^{m}}{m!}=\sum_{m=0}^{\infty}\frac{d^{m}}{dx^{m}}\sum_{n=0}^{\infty}\mathrm{Bel}_{n,\lambda}(t)\frac{x^{n}}{n!}\frac{y^{m}}{m!} \label{17}\\
&=\sum_{n=0}^{\infty}\frac{\mathrm{Bel}_{n,\lambda}(t)}{n!}\sum_{m=0}^{\infty}\Big(\frac{d}{dx}\Big)^{m}x^{n}\frac{y^{m}}{m!}=\sum_{n=0}^{\infty}\frac{\mathrm{Bel}_{n,\lambda}(t)}{n!}\sum_{m=0}^{n}\binom{n}{m}x^{n-m}y^{m}.\nonumber\\
&= \sum_{n=0}^{\infty}\frac{\mathrm{Bel}_{n,\lambda}(t)}{n!}(x+y)^{n}=e_{\lambda}\Big(t\big(e_{\lambda}(x+y)-1\big)\Big)\nonumber\\
&=e_{\lambda}\bigg(t\bigg(e_{\lambda}(x)e_{\lambda}\bigg(\frac{y}{1+\lambda x}\bigg)-1\bigg)\bigg). \nonumber
\end{align}
On the other hand, by \eqref{5}, we get 
\begin{align}
&e_{\lambda}\bigg(t\bigg(e_{\lambda}(x)e_{\lambda}\bigg(\frac{y}{1+\lambda x}\bigg)-1\bigg)\bigg)=\bigg(1+\lambda t\bigg(e_{\lambda}(x)e_{\lambda}\bigg(\frac{y}{1+\lambda x} \bigg)-1\bigg)\bigg)^{\frac{1}{\lambda}} \label{18}\\
&=\bigg(1+\lambda t\big(e_{\lambda}(x)-1\big)+\lambda t e_{\lambda}(x)\bigg(e_{\lambda}\bigg(\frac{y}{1+\lambda x}\bigg)-1\bigg)\bigg)^{\frac{1}{\lambda}}\nonumber\\
&=\Big(1+\lambda t(e_{\lambda}(x)-1)\Big)^{\frac{1}{\lambda}}\bigg(1+\lambda \frac{te_{\lambda}(x)\big(e_{\lambda}\big(\frac{y}{1+\lambda x}\big)-1\big)}{1+\lambda t(e_{\lambda}(x)-1)}\bigg)^{\frac{1}{\lambda}}\nonumber\\
&=e_{\lambda}\Big(t(e_{\lambda}(x)-1\big)\Big)e_{\lambda}\bigg(\frac{te_{\lambda}(x)\big(e_{\lambda}\big(\frac{y}{1+\lambda x}\big)-1\big)}{1+\lambda t(e_{\lambda}(x)-1)}\bigg)\nonumber\\
&=e_{\lambda}\Big(t(e_{\lambda}(x)-1\big)\Big)\sum_{k=0}^{\infty}\frac{(1)_{k,\lambda}}{k!} \frac{t^{k}e_{\lambda}^{k}(x)\big(e_{\lambda}\big(\frac{y}{1+\lambda x}\big)-1\big)^{k}}{\big(1+\lambda t(e_{\lambda}(x)-1)\big)^{k}}\nonumber \\
&=e_{\lambda}\Big(t(e_{\lambda}(x)-1\big)\Big)\sum_{k=0}^{\infty}(1)_{k,\lambda}\bigg(\frac{1}{1+\lambda t(e_{\lambda}(x)-1)}\bigg)^{k}t^{k}e_{\lambda}^{k}(x)\sum_{m=k}^{\infty}{m \brace k}_{\lambda}\frac{y^{m}e_{\lambda}^{-m\lambda}(x)}{m!}\nonumber\\
&=\sum_{m=0}^{\infty}\frac{y^{m}}{m!}\sum_{k=0}^{m}(1)_{k,\lambda}{m \brace k}_{\lambda}t^{k}\frac{e_{\lambda}^{k-m\lambda}(x)}{\big(1+\lambda t(e_{\lambda}(x)-1)\big)^{k}}e_{\lambda}\Big(t(e_{\lambda}(x)-1\big)\Big) \nonumber\\
&=\sum_{m=0}^{\infty}\frac{y^{m}}{m!}\sum_{k=0}^{m}(1)_{k,\lambda}t^{k}{m \brace k}_{\lambda}\sum_{j=0}^{\infty}F_{j,\lambda}^{(k)}(-\lambda t,k-m\lambda)\frac{x^{j}}{j!}\sum_{l=0}^{\infty}\mathrm{Bel}_{l,\lambda}(t)\frac{x^{l}}{l!}\nonumber\\
&=\sum_{m,n=0}^{\infty}\bigg(\sum_{k=0}^{m}\sum_{l=0}^{n}(1)_{k,\lambda}{m \brace k}_{\lambda}\binom{n}{l}t^{k}F_{n-l,\lambda}^{(k)}(-\lambda t,k-m\lambda)\mathrm{Bel}_{l,\lambda}(t)\bigg)\frac{x^{n}}{n!}\frac{y^{m}}{m!}. \nonumber
\end{align}
Therefore, by \eqref{17} and \eqref{18}, we obtain the following theorem. 
\begin{theorem}
For $n,m\ge 0$, we have 
\begin{displaymath}
\mathrm{Bel}_{n+m,\lambda}(t)= \sum_{k=0}^{m}\sum_{l=0}^{n}(1)_{k,\lambda}{m \brace k}_{\lambda}\binom{n}{l}t^{k}F_{n-l,\lambda}^{(k)}(-\lambda t,k-m\lambda)\mathrm{Bel}_{l,\lambda}(t).
\end{displaymath}
\end{theorem}
Letting $\lambda \rightarrow 0$ gives Spivey's recurrence relation for Bell polynomials (see \eqref{10}, \eqref{12}):
\begin{displaymath}
\phi_{n+m}(t)=\sum_{k=0}^{m}\sum_{l=0}^{n}{m \brace k}\binom{n}{l}t^{k}k^{n-l}\phi_{l}(t), 
\end{displaymath}
where $m,n$ are nonnegative integers. \par 
From \eqref{14}, we note that 
\begin{align}
&\sum_{m,n=0}^{\infty}F_{n+m,\lambda}(t)\frac{x^{n}}{n!}\frac{y^{m}}{m!}=\sum_{m=0}^{\infty}\Big(\frac{d}{dx}\Big)^{m}\sum_{n=0}^{\infty}F_{n,\lambda}(t)\frac{x^{n}}{n!}\frac{y^{m}}{m!} \label{19} \\
&=\sum_{n=0}^{\infty}\frac{F_{n,\lambda}(t)}{n!}\sum_{m=0}^{\infty}\Big(\frac{d}{dx}\Big)^{m}x^{n}\frac{y^{m}}{m!}=\sum_{n=0}^{\infty}\frac{F_{n,\lambda}(t)}{n!}\sum_{m=0}^{n}\binom{n}{m}x^{n-m}y^{m} \nonumber\\
&=\sum_{n=0}^{\infty}\frac{F_{n,\lambda}(t)}{n!}(x+y)^{n}=\frac{1}{1-t(e_{\lambda}(x+y)-1)}.\nonumber
\end{align}
On the other hand, by \eqref{14}, we get 
\begin{align}
& \frac{1}{1-t(e_{\lambda}(x+y)-1)}=\cfrac{1}{1-t\Big(e_{\lambda}(x)e_{\lambda}\big(\frac{y}{1+\lambda x}\big)-1\Big)}\label{20}\\
&=\cfrac{1}{1-t\big(e_{\lambda}(x)-1\big)-te_{\lambda}(x)\Big(e_{\lambda}\big(\frac{y}{1+\lambda x}\big)-1\Big)}\nonumber\\
&=\frac{1}{1-t(e_{\lambda}(x)-1)}\cfrac{1}{1-\cfrac{te_{\lambda}(x)\Big(e_{\lambda}\big(\cfrac{y}{1+\lambda x}\big)-1\Big)}{1-t(e_{\lambda}(x)-1)}}\nonumber\\
&=\frac{1}{1-t(e_{\lambda}(x)-1)}\sum_{k=0}^{\infty}t^{k}\bigg(\frac{1}{1-t(e_{\lambda}(x)-1)}\bigg)^{k}e_{\lambda}^{k}(x)\bigg(e_{\lambda}\bigg(\frac{y}{1+\lambda x}\bigg)-1\bigg)^{k} \nonumber \\
&=\frac{1}{1-t(e_{\lambda}(x)-1)}\sum_{k=0}^{\infty}\frac{t^{k}k!e_{\lambda}^{k}(x)}{(1-t(e_{\lambda}(x)-1))^{k}}\frac{1}{k!}\bigg(e_{\lambda}\bigg(\frac{y}{1+\lambda x}\bigg)-1\bigg)^{k}\nonumber\\ 
&=\frac{1}{1-t(e_{\lambda}(x)-1)}\sum_{k=0}^{\infty}\frac{t^{k}k!e_{\lambda}^{k}(x)}{(1-t(e_{\lambda}(x)-1))^{k}}\sum_{m=k}^{\infty}{m \brace k}_{\lambda}\frac{y^{m}e^{-m \lambda}(x)}{m!}\nonumber\\
&=\sum_{m=0}^{\infty}\frac{y^{m}}{m!}\sum_{k=0}^{m}t^{k}k!{m \brace k}_{\lambda}\bigg(\frac{1}{1-t(e_{\lambda}(x)-1)}\bigg)^{k}e_{\lambda}^{k-m\lambda}(x)\frac{1}{1-t(e_{\lambda}(x)-1)}\nonumber\\ 
&=\sum_{m=0}^{\infty}\frac{y^{m}}{m!}\sum_{k=0}^{m}t^{k}k!{m \brace k}_{\lambda}\sum_{j=0}^{\infty}F_{j,\lambda}^{(k)}(t,k-m\lambda)\frac{x^{j}}{j!}\sum_{l=0}^{\infty}F_{l,\lambda}(t)\frac{x^{l}}{l!}\nonumber\\
&=\sum_{m=0}^{\infty}\sum_{n=0}^{\infty}\frac{y^{m}}{m!}\frac{x^{n}}{n!}\bigg(\sum_{l=0}^{n}\sum_{k=0}^{m}t^{k}k!{m \brace k}_{\lambda}\binom{n}{l}F_{l,\lambda}(t)F_{n-l,\lambda}^{(k)}(t,k-m\lambda)\bigg).\nonumber 
\end{align}
Therefore, by \eqref{19} and \eqref{20}, we obtain the following theorem. 
\begin{theorem}
For $m,n\ge 0$, we have 
\begin{displaymath}
F_{n+m,\lambda}(t)=\sum_{k=0}^{m}\sum_{l=0}^{n}k!{m \brace k}_{\lambda}\binom{n}{l}t^{k}F_{n-l,\lambda}^{(k)}(t,k-m\lambda)F_{l,\lambda}(t). 
\end{displaymath}
\end{theorem}
Letting $\lambda \rightarrow 0$, we obtain the following identity (see \eqref{13}):
\begin{equation*}
F_{n+m}(t)=\sum_{k=0}^{m}\sum_{l=0}^{n}k!{m \brace k}\binom{n}{l}t^{k}F_{n-l}^{(k)}(t,k)F_{l}(t). 
\end{equation*} \par

\section{Conclusion} 
Recent research has extensively explored degenerate versions of special polynomials, numbers, and functions, including gamma functions and even umbral calculus. These explorations have utilized a wide range of tools, such as combinatorial methods, generating functions, $p$-adic analysis, probability theory, quantum mechanics, and operator theory. \par
In our study, we focused on degenerate versions of the Bell and Fubini polynomials, specifically the fully degenerate Bell polynomials $\mathrm{Bel}_{n,\lambda}(x)$ and the degenerate Fubini polynomials $F_{n,\lambda}(x)$. We successfully derived Spivey-type recurrence relations for these polynomials by employing generating functions. \par
Our future work will continue to investigate a variety of degenerate versions of special numbers, polynomials, and functions, with the aim of discovering their applications across physics, science, engineering, and mathematics.

\end{document}